\documentclass[reqno,10pt]{amsart}

\usepackage[english]{babel}
\usepackage{color}
\usepackage{enumitem}
\usepackage{amsmath}
\usepackage{amssymb}
\usepackage{amsfonts}
\usepackage{amsthm}

\newcommand{\Q}{\mathbb Q}
\newcommand{\N}{\mathbb{N}}
\newcommand{\R}{\mathbb{R}}

\newcommand{\Z}{\mathbb{Z}}
\newcommand{\E}{\mathbb{E}}

\newcommand{\f}{\rightarrow}

\newcommand{\euc}{\operatorname{euc}}

\newcommand{\codiam}{\operatorname{codiam}}
\newcommand{\diam}{\operatorname{diam}}

\newcommand{\stsys}{\operatorname{stsys}}
\newcommand{\sys}{\operatorname{sys}}

\newcommand{\Vol}{\operatorname{Vol}}

\newtheorem{thm}{Theorem}[section]

\newtheorem{margulis}[thm]{Abelian Margulis' Lemma}
\newtheorem{qbdt}[thm]{Quantitative Bounded Distance Theorem}
\newtheorem{prop}[thm]{Proposition}
\newtheorem{lem}[thm]{Lemma}

\theoremstyle{remark}
\newtheorem{exmp}[thm]{Example}
\newtheorem{rmk}[thm]{Remark}

\theoremstyle{definition}

\begin{document}

\title[QBD Theorem and Margulis' Lemma for $\mathbb{Z}^n$-actions]{Quantitative Bounded Distance Theorem \\and Margulis' Lemma for $\mathbb{Z}^n$-actions, \\with Applications to Homology}
\author{F. Cerocchi, A. Sambusetti}

\AtEndDocument{
\bigskip{\footnotesize 
F. Cerocchi, \textsc{Centro di Ricerca Matematica \lq\lq Ennio De Giorgi''} \par  
\textsc{Collegio Puteano, Scuola Normale Superiore Piazza dei Cavalieri, 3 I-56100 Pisa } \par
\textit{E-mail address}: \texttt{filippo.cerocchi@sns.it } \par
\addvspace{\medskipamount}
A. Sambusetti, \textsc{Dipartimento di Matematica ``G. Castelnuovo''}  \par  
\textsc{Sapienza Universit\`a di Roma, P.le Aldo Moro 5, I-00185 Roma } \par
\textit{E-mail address}: \texttt{sambuset@mat.uniroma1.it} 
}}

\maketitle

\vspace{-3mm}
{\bf Abstract.} We consider the stable norm 
associated to a discrete, torsionless abelian group of isometries $\Gamma \cong \Z^n$ of a geodesic space $(X,d)$. We show that the difference between the stable norm $\| \;\, \|_{st}$ and the distance $d$ is bounded by a constant only depending on the rank $n$ and on upper bounds for the diameter of $\bar X=\Gamma \backslash X$ and the asymptotic volume $\omega(\Gamma, d)$. 
We also prove that the upper bound on the asymptotic volume
is equivalent to a lower bound for the stable systole of the action of $\Gamma$ on $(X,d)$; for this, we  establish  a Lemma \`a la Margulis for $\mathbb{Z}^n$-actions, which gives optimal estimates of $\omega(\Gamma,d)$ in terms of  $\stsys(\Gamma,d)$, and vice versa, and characterize the cases of equality. 
Moreover, we   show that all the parameters $n, \diam(\bar X)$ and $\omega (\Gamma, d)$ (or $\stsys (\Gamma,d)$) are necessary to bound the difference $d -\| \;\, \|_{st}$, by providing explicit counterexamples for each case. \linebreak
As an application in Riemannian geometry, we prove that the number of connected components of any optimal, integral $1$-cycle in a closed Riemannian manifold $\bar X$ either is bounded by an explicit function of the first Betti number, $ \diam(\bar X)$ and $\omega(H_1(\bar X, \mathbb{Z}))$, or is  a sublinear function of the mass.


\section{Introduction}

Consider a geodesic metric space $(X,d)$ with a {\em $\Z^n$-periodic metric}, i.e. admitting a discrete, torsionless abelian group of isometries $\Gamma$ 
of rank $n$
acting 
properly discontinuously:
we mainly think to the Cayley graph of a word metric on $\Z^n$, or to a $\Z^n$-covering of a compact Riemannian or Finsler manifold.
A motivating example is the torsion free homology
covering $X$ of any compact manifold $\bar X$ with nontrivial first Betti number, which has  automorphism group $\Gamma=H_1(\bar X,\Z)/tor$. 





\noindent The associated {\em stable norm} on $\Gamma$ is defined as:
\vspace{-3mm}

$$ \| \, \gamma \, \|_{st} = \lim_{k\rightarrow \infty} \frac{1}{k} d (x_0, \gamma^k. x_0)$$

\noindent and clearly it does not depend on the choice of $x_0 \in X$, by the triangular inequality. An isomorphism $\Gamma \cong \Z^n$ being chosen, this yields a well-defined norm\footnote{since it can be bounded from below by a multiple of a word metric, cp. \cite{BBI}.} on $\R^n$,
extending the definition by homogeneity to $\Q^n$ first, and then to real coefficients by uniform continuity. 
For instance, when $\Gamma=H_1(\bar X,\Z)/tor$ is the automorphism group of  the torsion free homology covering of a compact Riemannian manifold $\bar X$, the stable norm of coincides with the norm induced by the Riemannian length in the homology with real coefficients,   that is:
\vspace{-3mm}

\small
$$  \| \, \gamma \,\|_{st} 
= \inf \left\{ \sum_k \! |a_k| \ell (\gamma_k) \; : \, a_k \!\in\! \R, \gamma_k \mbox{ Lipschitz $1$-cycles, } \gamma \!=\!\sum_k a_k \gamma_k \mbox{ in } H_1(\bar X,\R) \right\} $$
\normalsize

\noindent see \cite{federer}, \cite{gro2}.
\vspace{2mm}

It is folklore ({\em Bounded Distance Theorem}, cp. \cite{buriv2},\cite{BBI}, \cite{gro2}) that, when $\Gamma$ acts {\em cocompactly},  
then   $(X,d)$ is almost isometric to $(\Z^n,  \| \|_{st})$:
namely, for every $x_0 \in X$, there exists a constant $C$ such that 
$$\left| \, d (x_0,\gamma . x_0) - \|\,\gamma\,\|_{st} \, \right| < C \mbox{ for all } \gamma \in \Gamma.$$
This fact was originally  proved D. Burago for periodic metrics on $\R^n$ (see \cite{Bur}, \cite{gro1});
however, we were not able to find a complete proof of the general case in literature.
The first purpose of this note is to  investigate to what extent the constant $C$ depends on the basic geometric invariants of $X$, 
i.e. to estimate how far a space admitting an abelian action is from a normed vector space.
We prove:

\begin{qbdt}
\label{QBDT}
Let $\Gamma = \Z^n$ act freely and properly discontinuously by isometries on a length space  $(X,d)$, 
with compact quotient. There exists a constant $c=c(n,D,\Omega)$ such that for all $x_0 \in X$
\vspace{-3mm}

\begin{equation}
\label{QBDT-eq}
\left| \, d(x_0,\gamma x_0)-\|\,\gamma \,\|_{st} \, \right|< c(n,D,\Omega)
\end{equation}

\noindent where $D$ and $\Omega$ are, respectively, upper bounds for the codiameter 
and   the asymptotic volume of $\Gamma$ with respect to $d$.
\end{qbdt}

\noindent We call {\em co-diameter} of  $\Gamma$
the diameter of the quotient  \nolinebreak  $\bar X =\Gamma \backslash X$. \\
The {\em asymptotic volume} of a group $\Gamma \cong \Z^n$, endowed with a $\Gamma$-invariant metric $d$, is the asymptotic invariant defined as 
\vspace{-3mm}

$$ \omega (\Gamma, d) 
= \lim_{R\rightarrow \infty} \frac{ \sharp \{ \gamma \, : d(x,\gamma . x) < R\}}{R^n} $$ 


\noindent For $\Gamma$ acting on $(X,d)$ as above, any choice of a base point $x_0 \in X$ yields  a left-invariant distance $d_{x_0}$ on $\Gamma$, by identification with the orbit $\Gamma .x_0$; clearly, the asymptotic volume $\omega(\Gamma,d_{x_0})$ does not depend on $x_0$, and we shall simply write $\omega(\Gamma,d)$. 
Moreover, given any  $\Gamma$-invariant measure $\mu$ on $X$,  it is easy to see, by a packing argument, that it equals the usual asymptotic volume of the measure metric space $(X,d, \mu)$  divided by the measure of the quotient, i.e.:
$$\omega_\mu (X,d) = \lim_{R \rightarrow \infty }  \frac{ \mu \left( B_{(X,d)} (x_0, R) \right) }{R^n}    =  \mu (\bar X) \cdot  \omega (\Gamma, d).$$

\noindent As a consequence of  the QBD Theorem \ref{QBDT}, we have an explicit control of the growth function of balls and annuli in $(\Gamma, d)$ (cp. Proposition \ref{v-estimates}), and
of the Gromov-Hausdorff distance between $(X,\lambda d)$ and its asymptotic cone $( \R^n, \| \;\;\|_{st})$ in terms of $n,D,\Omega$:
\vspace{-3mm}

$$d_{GH} \bigg( (X,\lambda d), \;(\R^n,\|\;\|_{st}) \bigg) \leq \lambda \cdot (c+2D) $$ 

\noindent (notice that, for abelian groups endowed with a word metric, the linearity of the rate of convergence of $(X,\lambda d)$ to the asymptotic cone was already known, cp. \cite{buiv}). \\



The QBD Theorem \ref{QBDT} is obtained combining Burago's original idea with a careful control of the dilatation of  ``natural'' maps
$ (\R^n, euc)  \ \leftrightarrows  (X,d)$ quasi-inverse one to each other\footnote{This difficulty does not emerge in \cite{Bur}, where $\Gamma \curvearrowright \R^n$, since in that case we have two metrics on a torus, which are obviuosly bi-Lipschitz via the identity map.}.
More precisely, the maps  are induced from the identification of $\Z^n$ with a finite index subgroup ${\mathcal Z}$ of $\Gamma$ generated by a set $\Sigma_n$ of $n$ linearly independent vectors $(\gamma_k)$. 
The bounds on the codiameter and the asymptotic volume are then needed to control the index
 $[ \Gamma : {\mathcal Z}]$ and the relative variation of $d/d_{\Sigma_n}$ on ${\mathcal Z}$. \linebreak
  For this, we prove in Section \S3 a Lemma \`a la Margulis\footnote{The classical Margulis' Lemma, in negative curvature, gives an estimate (at some point $x_0$) of the minimal displacement of a group $\Gamma$ acting on a Cartan-Hadamard manifold $X$, under a lower bound on the curvature of $X$. It has been extended in several directions, in particular with a bound on the volume entropy replacing the bound on curvature, cp. \cite{bcg}, \cite{cerocchi}.}
for abelian groups, which gives an estimate from below of the minimal displacement of $\Gamma$ in terms of an upper bound on the asymptotic volume. 
 Namely, let $\Gamma^\ast = \Gamma \setminus \{ e\}$ and define, respectively, the  systole and  the stable systole of the action of $\Gamma$ on $(X,d)$ as
 \vspace{-2mm}
 
 $$\sys (\Gamma, d) = \inf_{x\in X} \inf_{\gamma\in\Gamma^\ast} d(x,\gamma . x)$$
 $$\stsys (\Gamma, d)=\inf_{\gamma \in \Gamma^\ast} \| \,\gamma \,\|_{st}$$
 
\noindent Then clearly $\stsys (\Gamma,d) \leq \sys (\Gamma,d)$, and we prove:

\begin{margulis}\label{margulis}
Let $\Gamma = \Z^n$ act  freely and properly discontinuously by isometries on a length space $(X,d)$, with cocompact quotient. Then:
\vspace{-3mm}

\begin{equation}
\label{eqmargulis}
\frac{2}{n!}   \cdot \;\;\frac{1}{\codiam (\Gamma,d)^{n-1}   \cdot  \omega (\Gamma, d)}
\le \stsys (\Gamma,d) \le
 \frac{2}{\omega (\Gamma, d)^{1/n}}
\end{equation}

\vspace{1mm}
\noindent Moreover, these inequalities are optimal and the equality cases characterize,  up to  almost-isometric equivalence, the action of specific lattices of $\R^n$, endowed with particular polyhedral norms. 
Namely, if $\codiam(\Gamma,d) = D$, $\stsys (\Gamma, d)=\sigma$ and $\omega ( \Gamma, d) \leq \Omega$,  then there exists a constant $C=C(n,D,\Omega)$ such that:
\vspace{1mm}

\noindent $\bullet$  the equality holds  in the left-hand side if and only if there is an equivariant, \linebreak $C$-almost isometry 
 $f : (X,d)  \rightarrow  (\R^n, \| \; \|_{1})$, with respect to the action by translations of the lattice 
 $\Gamma_0 = \sigma \!\cdot\! \Z \times 2D \!\cdot \! \Z^{n-1} \cong \Gamma$ on $\R^n$;
\vspace{1mm}

\noindent $\bullet$ the  equality holds  in the right-hand side if and only if there is an equivariant, \linebreak $C$-almost isometry $f:  (X,d) \! \rightarrow \!   (\R^n\!, \| \; \|_{h}  )$, 
   with respect to the canonical action of $\Gamma$ on $\R^n$
 and where $\| \; \|_{h}$ is a  
  {\em parallelohedral norm}\footnote{Examples of parallelohedral norms are:\\
\noindent (i) in dimension $n=2$, all norms whose unit ball is either a parallelogram or a convex hexagon with congruent  opposite sides
(these are the only convex polygons which tasselate $\R^2$ under the action by translations of a $2$-dimensional lattice, cp. \cite{fedorov1885}, \cite{shephard78}, \cite{shephard80}); \\
\noindent (ii) in dimension $n=3$, there are precisely 37 types of parallelohedra, cp. \cite{fedorov1899}, including for instance the standard $n$-cube (which gives rise to the  sup-norm) or those obtained from $\Z^n$-tessellation by  prisms with $2$-dimensional base as in (i). \\
 \noindent The complete classification of parallelohedral norms is a particular case of Hilbert's eighteenth problem (tiling the Euclidean space by congruent polyhedra) and will be not pursued further here.}; 
 that is, a norm whose unit ball is a {\em $\Gamma$-parallelohedron} (a convex polyhedron  which tiles $\R^n$ under the  action by translations of $\Gamma$, i.e. whose $\Gamma$-translates cover $\R^n$ and have disjoint interiors).
\vspace{1mm}

\end{margulis}

 The left-hand side of (\ref{eqmargulis}) shows that,  provided that the co-diamenter is bounded,  an upper bound of the  asymptotic volume is equivalent to a lower bound of the stable systole. Therefore, the constant $c(n,D, \Omega)$ in Theorem \ref{QBDT} can as well be expressed in terms of rank, co-diameter and of a lower bound $\stsys (\Gamma,d) \geq \sigma$, instead of an upper bound $\omega(\Gamma,d) \leq \Omega$.\\
Notice that $\stsys (\Gamma,d)$  cannot be bounded below uniquely in terms of $n$ and $\omega (\Gamma,d)$: any flat Riemannian torus $(T,euc)$ with unitary volume has fundamental group $\Gamma=\pi_1(T)$ with   asymptotic volume equal to 
the volume $\omega_n$  of the unit ball in $\E^n$, but the systole of $\Gamma$  can be arbitrarily small (provided that the diameter of $T$ is sufficiently large).

\pagebreak

\vspace{2mm}
It is natural to ask whether one can drop the dependence of the constant $c$ in the QBD Theorem on any of the parameters $n,D, \Omega$ or $\sigma$, and possibly replace the dependence on the stable systole by a lower bound on the systole. In Section \S4 we give counterexamples ruling out each of these possibilities. In particular, one cannot generally bound $d-\| \;\, \|_{st}$ only in terms of rank, co-diameter and systole: there exists a sequence of actions of $\Z^n$ on length spaces $(X_k,d_k)$ with $\sys(\Z^n, d_k)\geq 1$ and $\codiam (\Z^n, d_k) \leq 1$ such that the difference between $d_k$ and the corresponding stable norms $|\, \; |_{st,k}$ is arbitrarily large, cp. Example \ref{exdiameter}. 
The same example also shows that a lower bound of the systole does not imply any upper bound for the asymptotic volume, i.e. the right-hand side of (\ref{eqmargulis}) does not hold with the stable systole replaced by the systole.


\vspace{2mm}
Finally, in section \S\ref{sectioncc} we use the QBD Theorem to address the following basic problem  on a closed Riemannian manifold $\bar X$:  given an  integral homology class $\gamma \in  H_1 (\bar X, \Z)$,   what is the minimal number $\#_{CC}$ of connected components of an   optimal  cycle in $\gamma$? Namely, we want to estimate the number 
\vspace{-2mm}

$$N(\gamma) = \min \{ \#_{CC} (  \textsf{c}  ) \; | \;  \textsf{c}  \in Z_1 (\bar X, \Z),  \;  [ \textsf{c} ] =\gamma, \; \ell( \textsf{c} ) = |\gamma |_{H_1} \}$$

\vspace{1mm}
\noindent where $|\gamma |_{H_1}$ is the {\em mass}  in homology, i.e. the total length  of a 
shortest\footnote{Notice that a $1$-cycle of minimal length in $\bar X$ always exists, and is given by a finite collection of closed geodesics, by general representation results of minimizers in homology by currents, and by regularity  of rectifiable $1$-currents.},  
possibly disconnected, collection of closed curves  representing $\gamma$.
We call  {\em optimal}  a cycle $\textsf{c} \in [\gamma ]$  which is length-minimizing  in its class and having  precisely the minimum number $N(\gamma)$ of connected components. 


\noindent Recall that the {\em homological systole} $\sys_{H_1} (\bar X) $ of $\bar X$ is the length of the shortest closed geodesic which is non-trivial in homology; if  $\Gamma={H_1} (\bar X, \Z) $   and $(X,d)$ is the Riemannian homology covering of $\bar X$,   we clearly have $\sys_{H_1} (\bar X)  = \sys  (\Gamma, d)$. 
Notice that  a lower bound of the homological systole 
$\sys_{H_1} (\bar X) \geq \sigma_1$  
(as given for instance, in the torsionless case, by the left-hand side of the Abelian Margulis' Lemma \ref{margulis}) readily implies  an estimate $N(\gamma) \leq \sigma_1^{-1} |\gamma |_{H_1}$.
However,  as an application of  the QBD Theorem, we actually   show that $N(\gamma) $ is   sublinear  in  $|\gamma |_{H_1}$:

\begin{thm}
\label{teorcc}
Assume that $\bar X$ has  first Betti number $b_1(\bar X) =n$,  $\diam(\bar X) < D$ and $\omega (H_1(\bar X, \Z)) < \Omega$.
Then, for any  torsionless  homology  class $\gamma \in H_1  (\bar X, \Z)$: 
\vspace{1mm}

\noindent (i) either $N( \gamma)$ is bounded by an explicit, universal function $N(n,D,\Omega)$,
\vspace{1mm}

\noindent (ii) or $\displaystyle N( \gamma) \leq 2 \cdot 3^{2n} \cdot \Omega^{\frac{1}{n+1}} \cdot |\gamma |_{H_1
}^{\frac{n}{n+1}}$.
\end{thm}

\vspace{-1mm}
\noindent We will see that one can take $N(n,D,\Omega) = 2^{18 n^3} \! \cdot \! n^{2n} \! \cdot \!  (n!)^{n(n+2)} \!  \cdot (\Omega D^n+1)^{6n^2}$.

\noindent It is noticeable that the bound (ii) does not even depend on the diameter of $\bar X$.

\vspace{2mm}
\small {\sc Aknowledgements}. We thank D. Massart and S. Saboureau for   useful discussions.
\normalsize

\section{QBD Theorem.}

Let $\Gamma\cong\Z^n$ act freely, properly discontinuously and cocompactly by isometries on a length space $(X,d)$. For any given $x_0\in X$ we consider the left invariant metric  on $\Gamma$ given by $d_{x_0}(\gamma_1\,,\,\gamma_2)= d(\gamma_1 . x_0\,,\,\gamma_2.x_0)$.  
 We will write  $d_S$ for the word metric relative to a generating set $S$ of a group, and also use the abridged notations  
 $| \gamma |_{x_0}= d(x_0,\gamma.x_0)$, $| \gamma |_S= d_S(e,\gamma)$. \pagebreak 

 Assume that $\diam(\Gamma\backslash X)\le D$, $\omega(\Gamma, d)\le\Omega$ and $\sys(\Gamma,d)\ge\sigma$.
 
\vspace{1mm}
\noindent We consider the generating 
set\footnote{The elements with $d(\gamma x_0, x_0) \leq 2D$ suffice to generate $\Gamma$, 
cp. \cite{gro2}, p.91; the constant $3D$ is chosen here to bound from below $d_{x_0} / d_\Sigma$.} 
$\Sigma_D=\{\gamma\in\Gamma^*\;|\;  d( \gamma x_0, x_0) \le 3D\}$, and we extract from $\Sigma_D$ a set $\Sigma_n=\{\gamma_1,...,\gamma_n\}$ of $n$ linearly independent vectors which generate a finite index subgroup $\mathcal Z=\langle\Sigma_n\rangle$, again isomorphic to $\Z^n$. Then,   fix once and for all  a set of representatives $S=\{s_0=e,s_1,..., s_d\}$  for $\Gamma/\mathcal Z$ which are \textit{minimal} for the word metric $d_{\Sigma_D}$ associated to the generating set $\Sigma_D$ of $\Gamma$. 

\vspace{1mm}
\noindent  Let us consider the map $f:(\mathcal Z, d_{x_0})\f(\Z^n, \euc)$ defined by sending each $\gamma_i$ to the \textit{i-th} vector of the standard basis of $\R^n$.
We shall prove that $f$ and $f^{-1}$ are two Lipschitz maps, whose Lipschitz constants  $M$ and $M'$   are bounded in terms of our geometric data $n,D, \Omega$ and $\sigma$; we will then  extend  $f$  to a Lipschitz map $F:(X, d)\f (\R^n,\euc)$.
The purpose of the next lemmas is to estimate the constants $M$, $M'$ by comparing with the   dilatations of the following maps
\begin{equation}\label{maps}
f: (\mathcal Z, d_{x_0})\f (\mathcal Z, d_{\Sigma_D}|_{\mathcal Z})\f (\mathcal Z, d_{\widehat{\Sigma}_n})\f (\mathcal Z, d_{\Sigma_n})\f(\Z^n, \euc)
\end{equation}
where:
 
\begin{itemize}[leftmargin=*]
 \item  $d_{\Sigma_D}|_{\mathcal Z}$ is the restriction of the word metric $d_{\Sigma_D}$ to $\mathcal Z$\,;
\item  $d_{\Sigma_n}$ is the word metric on $\mathcal Z$ relative to $\Sigma_n$\,;
\item  $d_{\widehat{\Sigma}_n}$ is the word metric on $\mathcal Z$ relative to the generating set $\widehat{\Sigma}_n$ of ${\mathcal Z}$ defined by
$$\widehat{\Sigma}_n=\{ s_i\,\sigma\, s_j^{-1} \;|\; s_i \in S,\, \sigma\in\Sigma_D \mbox{ and }   s_i\,\sigma\, s_j^{-1} \in \mathcal Z^* \}\,.$$
\end{itemize}

\noindent (with \small $B_{\! (\mathcal Z, d_{x_0})}\! (r), B_{\! (\mathcal Z, \Sigma_D)}\! (r), B_{\! (\mathcal Z, \Sigma_n)}\! (r)$  \normalsize  and \small $B_{\! (\mathcal Z, \widehat{\Sigma}_n)}\! (r)$ \normalsize the relative balls  centered at $e$.)

\noindent Notice that $\Sigma_n \subset \widehat{\Sigma}_n$ (since $s_0=e \in S$), 
but we might have  $\Sigma_D \not\subset \widehat{\Sigma}_n$. \\
Moreover, remark that $(\mathcal Z, d_{\Sigma_n})$ is isometric to $\Z^n$ endowed with the canonical word metric $\| \; \|_1$, so we have:
 \vspace{-2mm}
 
\begin{equation}
 \label{eqcanonicaldistance}
 \frac{1}{\sqrt{n}} \cdot |\,\gamma\,|_{\Sigma_n}\le \| f(\gamma) \|_{euc} \le  |\,\gamma\,|_{\Sigma_n}
\end{equation}
 
 \vspace{-2mm}
 
\begin{equation}
 \label{eqcanonicalvolume}
 \omega (\mathcal Z, d_{\Sigma_n}) = \frac{2^n}{n!}
\end{equation}
 
\vspace{2mm}
\begin{lem}[cp. \cite{gro2}]\label{gromov}
The set $\Sigma_D$ is a generating set for $\Gamma$ such that:
$$\frac{\sigma}{2}\cdot |\,\gamma\,|_{\Sigma_D}\le d(x_0,\gamma.x_0)\le 3D\cdot |\,\gamma\,|_{\Sigma_D}$$
\end{lem}

\begin{lem}\label{geo-alg}
For all $\gamma \in \mathcal Z$ we have:
$$| \gamma |_{\widehat{\Sigma}_n}\le |\gamma|_ {\Sigma_D}  
\le(2\,[\Gamma :\mathcal Z]+1)\cdot |\gamma |_{\widehat{\Sigma}_n}$$
\end{lem}

\textbf{Proof.} Let $\gamma=\gamma_1\cdots\gamma_\ell \in \mathcal Z $ with $\gamma_i\in\Sigma_D$. Assume that $\gamma_1\cdots\gamma_i\,\mathcal Z= s_{k_i}\,\mathcal Z$, then $\gamma$ can be written as:
$$\gamma=(s_{k_0}\gamma_1\cdot s_{k_1}^{-1})\cdot(s_{k_1}\,\gamma_2\,s_{k_2}^{-1})\cdots (s_{k_{\ell-2}}\,\gamma_{\ell-1}\,s_{k_{\ell-1}}^{-1})\cdot (s_{k_\ell -1}\,\gamma_\ell s_{k_l})$$
with $s_{k_0}\!=\!s_{k_l}\!=\!e$, and any $s_{k_{i-1}}^{-1}  \gamma_i  s_{k_i}$ either is trivial or belongs to  $\widehat{\Sigma}_n$. \\ Therefore $ |\,\gamma\,|_{\widehat{\Sigma}_n} \le |\,\gamma\,|_{\Sigma_D}$. 
For the second inequality,  recall that any class $s_i\,\mathcal Z$ can be written as $s_i\,\mathcal Z=\gamma_1\cdots\gamma_k\,\mathcal Z$ with $\gamma_i\in\Sigma_D$ and  $k\le [\Gamma:\mathcal Z]$. So,  every representative $s_i$, being $\Sigma_D$-minimal, satisfies $|\,s_i\,|_{\Sigma_D}\le [\Gamma : \mathcal Z]$, which implies   $|\,\gamma\,|_{\Sigma_D}\le (2\,[\Gamma :\mathcal Z]+1)\cdot |\,\gamma\,|_{\widehat{\Sigma}_n}$. $\Box$\\

\begin{lem}\label{est_quot} 
The subgroup $\mathcal Z$ satisfies:
\vspace{1mm}

\noindent (i) $\;\;[\Gamma :\mathcal Z]\le\frac{n!}{2^n}\Omega  (3D)^n $\,;
\vspace{1mm}

\noindent (ii) $\;d_{x_0}(\gamma, \mathcal Z)\le  \frac{n!}{2^n}\Omega (3D)^{n+1} $, for any $\gamma\in\Gamma$\,;
\vspace{1mm}

\noindent (iii) $\diam( \mathcal Z \backslash X) \leq D+    \Omega\frac{n!}{2^n}   (3D)^{n+1}$.

\end{lem}
\vspace{1mm}
\textbf{Proof.} We consider the set $S=\{s_i\}_{i=0,...d}$ of  representatives of $\Gamma/\mathcal Z$ with minimal $\Sigma_D$-length. Let $M=\max_{s_i} |\,s_i\,|_{\Sigma_D}$. Then
$$\# B_{(\Gamma,  {\Sigma_D})}(R)\ge [\Gamma:\mathcal Z]\cdot \# B_{(\mathcal Z,  {\Sigma_D})}(R-M)\ge [\Gamma : \mathcal Z]\cdot \# B_{(\mathcal Z, d_{\Sigma_n})}(R-M)\;.$$
Dividing by $R^n$ and   taking the limit for $R\f+\infty$ yields 
$[\Gamma : \mathcal Z]\le\frac{\omega(\Gamma, d_{\Sigma_D} )}{\omega(\mathcal Z, d_{\Sigma_n})}\;.$
\linebreak
By  Lemma \ref{gromov} we have $\omega(\Gamma, d_{\Sigma_D})\le (3D)^n \omega(\Gamma, d)  \leq (3D)^n \Omega$, while $\omega(\mathcal Z, d_{\Sigma_n})=\frac{2^n}{n!}$  by  (\ref{eqcanonicalvolume}); this proves  (i).\\
To prove (ii), notice that the set $\{ \gamma   \mathcal Z \; | \; \gamma  \in \Sigma_D\} $  generates $\Gamma/\mathcal Z$, and that every class $s_i  \mathcal Z$ is   product of at most $[\Gamma : \mathcal Z]$ classes $\gamma_i  \mathcal Z$ with $\gamma_i \in\Sigma_D$. Since any element of $\gamma\in \Gamma$ lies in some  coset $s_i \mathcal Z$, the $\Sigma_D$-distance of $\gamma$ from $\mathcal Z$ is at most  $ [\Gamma :\mathcal Z]$. Then,  Lemma \ref{gromov} yields $d_{x_0}(\gamma, \mathcal Z)\le 3D\cdot [\Gamma : \mathcal Z]\,$. \\
Assertion (iii) then follows from (ii), as $\diam (\Gamma \backslash X) \le D .\quad \Box$
\vspace{2mm}

\begin{lem}\label{alg-alg}
The generating set $\widehat{\Sigma}_n$ of $\mathcal Z$ satisfies: 
\vspace{1mm}

\noindent (i)  $\;\; \omega ( \mathcal Z, d_{ \widehat{\Sigma}_n  }) 
\le (2^{n+3} n! )^n \cdot \Omega D^n \cdot (\Omega D^n+1) \, $;
\vspace{1mm}

\noindent (ii)  $\; |\hat\gamma |_{\Sigma_n}\le L(n,D,\Omega)=2^{n^2+4n+3} (n!)^{n+1} \cdot \Omega D^n \cdot  (\Omega D^n+1)^n \;$ for all $\hat\gamma\in\widehat{\Sigma}_n$;
\vspace{2mm}

\noindent (iii)  $| \gamma  |_{\Sigma_n} \le L(n,D,\Omega)\cdot | \gamma |_{\widehat{\Sigma}_n}$
for all $\gamma\in \mathcal Z$.
 \end{lem}

\vspace{1mm}
\textbf{Proof.} 
 By the Lemmas \ref{gromov} and \ref{geo-alg} we have 
 $| \gamma |_{\widehat{\Sigma}_n} \geq \frac{1}{ 3D  (2\,[\Gamma :\mathcal Z]+1)} d(\gamma x_0, x_0)$, hence
 $   \omega ( \mathcal Z, d_{ \widehat{\Sigma}_n  })  
 \le \left[ 3D  (2\,[\Gamma :\mathcal Z]+1) \right]^n \omega  (\Gamma, d_{x_0})$, so (i) follows from Lemma \ref{est_quot}.  \\
To prove (ii),  assume that $\hat\gamma \in\widehat{\Sigma}_n$ has $\Sigma_n$-length $\ell$, 
so it can be written as a product \linebreak $\hat\gamma=\gamma_{i_1}\cdots\gamma_{i_\ell}$, with every $\gamma_{i_k} \in \Sigma_n$.
The sequence $(\gamma_{i_1},...,\gamma_{i_\ell})$ corresponds to a geodesic path $c_0$ in the Cayley graph $\mathcal C(\mathcal Z,\Sigma_n)$.   Let $c$  be the path in  $\mathcal C(\mathcal Z ,\Sigma_n)$  obtained by concatenation of all the paths $c_k=\hat \gamma^k . c_0$; notice that, since $(\mathcal Z, \Sigma_n)$ is isometric to $(\Z^n, | \; |_1)$, the path $c$ is still geodesic. 
Consider now a new generating set: $\Sigma_{n}({\hat \gamma})=\Sigma_n\cup\{\hat\gamma\} \subset \hat \Sigma_n$ and call for short $d_{\hat \gamma}$   the corresponding word metric. 
Chosen a radius $R = m \ell$, for $m > 0$, we consider the points  $P_i = \hat \gamma^{2mi}$ on the geodesic $c$, and we remark that:
\vspace{-3mm}

$$ \bigsqcup_{i=0}^{\lfloor \ell/2 \rfloor} B_{ (\mathcal Z, d_{\Sigma_n)}} (P_i, R-2mi) 
\subset B_{(\mathcal Z, d_{\hat \gamma})} (e,R)$$
Actually, for  $j \neq i \leq \ell/2$ the balls  $B_{ (\mathcal Z, d_{\Sigma_n})} (P_j, R-2mj)$ and  $B_{ (\mathcal Z, d_{\Sigma_n)}} (P_i, R-2mi)$ are disjoint, since 
$d_{\Sigma_n}(P_i, P_j) \geq  |  \hat\gamma^{2m} |_{\Sigma_n} = 2m \ell =2R$; moreover, these balls are all contained in  $B_{(\mathcal Z, d_{\hat \gamma} )} (e,R)$ as 
  $d_{\hat \gamma}  (e, P_i) = | \hat \gamma^{2mi} |_{\Sigma_{n}({\hat \gamma})} \leq 2mi$. \\
  Also, notice that, as 
  $\omega (\mathcal Z, d_{\Sigma_n}) = \frac{2^n}{n!}$, we have 
  \vspace{-3mm}
  
 $$\# B_{ (\mathcal Z, \Sigma_n)} (P_j, R-2mi) 
 = \# B_{ (\mathcal Z, \Sigma_n)} (m(\ell-2i))  \geq  \frac{2^{n-1}}{n!} m^n (\ell -2i)^n$$
  for $m \gg 0$.  Thus:
 \vspace{-4mm} \small
  
$$\# B_{ (\mathcal Z, d_{\hat \gamma}) }(R)
\ge \sum_{i=0}^{\lfloor \ell/2 \rfloor} \frac{2^{n-1}}{n!} m^n (\ell -2i)^n
\ge \frac{2^{n-1}}{n!} m^n \ell^n \sum_{i=0}^{\lfloor \ell/3 \rfloor} \left(\frac23 \right)^n
\ge \left(\frac43 \right)^n \frac{\ell R^n}{ 6 \cdot n!}
$$

\normalsize \noindent
which shows that   $\omega ( \mathcal Z, d_{\hat \gamma})  \geq  \left(\frac43 \right)^n \frac{\ell}{ 6 \cdot n!}$.
On the other hand, we know by the above lemmas 
that 
$ | \;\; |_{\Sigma_n({\hat \gamma})} \ge | \;\; |_{\hat{\Sigma}_n} 
\geq \frac{2^{n-1}}{3D(n! \cdot \Omega(3D)^{^{n}} + 2^{^n})} \cdot | \;\; |_{x_0}$,
so
$$\left(\frac43 \right)^n \frac{\ell}{ 6 \cdot n!} \le \omega ( \mathcal Z, d_{\hat \gamma} )
\le (2^{n+3} \cdot n!)^n \cdot \Omega D^n(\Omega D^n+1)^n$$ 
which gives (ii).
 The third statement clearly follows from  (ii).$\Box$

\vspace{5mm}
We deduce by  the lemmas above
that the map $f$ defined in (\ref{maps}) is a  bi-Lipschitz map, with Lipschitz constants  given by:
\small
\begin{equation}
\label{eqM} 
M(f)\le M=M(n,\Omega, D,\sigma)=\frac{1}{\sigma}\cdot  2^{n^2+4n+4}\cdot (n!)^{n+1} \cdot  \Omega D^n (\Omega D^n+1)^n
\end{equation}

\vspace{-4mm}

\begin{equation}
\label{eqM'} 
M(f^{-1})\le M'=M'(n,\Omega, D)=  8\, \left( \frac32 \right)^n \sqrt{n} \cdot n!  \cdot \Omega D^{n+1}
\end{equation}
\normalsize
We will prove in the next section that we can get rid of the dependence on $\sigma$.\\
 Now, we extend $f$ to a $M''$-Lipschitz map $F:(X,d)\f (\R^n, euc)$, with $M''=\sqrt{n}M$, by    extending each coordinate function $f_i$ of $f$ as follows
$$F_i(x)=\inf_{\gamma\in\mathcal Z} \left( f_i(\gamma.x_0)+ M(f)\cdot d(\gamma.x_0, x) \right) $$

\vspace{-2mm}
\noindent (notice that each $F_i$ is $M$-Lipschitz, and then $F$ is $\sqrt{n} M$-Lipschitz.)\\

\subsection{End of the proof of  Theorem \ref{QBDT}.} 
\label{sectionproof}
${}$

\noindent We switch now to the additive notation for the abelian groups $\Gamma$ and $\mathcal Z$, for easier comparison with $\Z^n$.  Assume first that $\gamma\in\mathcal Z$, and let $ \textsf{c}: I=[0,\ell]\f X$ be a minimizing geodesic (i.e. $d( \textsf{c}(t),  \textsf{c}(t'))= |t-t'|$) from $x_0$ to $2\gamma .x_0$. Then, we apply the following lemma due to  D. Burago and G. Perelman to the path $ \textsf{c}_o=F\circ  \textsf{c}$, going from the origin $o$ of $\R^n$ to $2f(\gamma) . o$ \;:

\begin{lem}[D. Burago, G. Perelman]
Let $ \textsf{c}:I=[0,\ell]\f\R^n$ be a Lipschitz path. There exists an open set $A=\bigcup_{i=1}^m (a_i,b_i)\subset [0,\ell]$ with $m\le n$ and with Lebesgue measure $\lambda(A)\le\frac12 \lambda(I)=\frac{\ell}{2}$ such that
\vspace{-2mm}
\small

$$\sum_{i=1}^m ( \textsf{c}_o(b_i)- \textsf{c}_o(a_i))= \frac{ \textsf{c}_o(\ell) - \textsf{c}_o (0)}{2}$$
\normalsize
\end{lem}

\noindent This Lemma provides a new path $\frac{ \textsf{c}_o}{2}: J=[0,\lambda(A)]\f\R^n$ going from the origin \nolinebreak $o$ to $f(\gamma).\,o=\sum_1^m ( \textsf{c}_o(b_i)- \textsf{c}_o(a_i))$, defined   concatenating   the paths $ \textsf{c}_{o, i}= \textsf{c}_o|_{[a_i, b_i]}$
$$\frac{ \textsf{c}_o}{2}= \textsf{c}_{o,1} *\cdots * \textsf{c}_{o,m} \; -  \textsf{c}_{o}(a_1)$$
(where $\alpha*\beta$ in $\R^n$ means that the path $\beta$ is translated in order that its origin coincides with the endpoint of $\alpha$). 
 Consider now, for each $i=1,..,m$, the orbit points $\alpha_i. x_0$, $\beta_i. x_0\in\mathcal Z x_0$ closest respectively to $ \textsf{c}(a_i),\, \textsf{c}(b_i)\,\in X$, and let $ \textsf{c}_i$ be a  minimizing geodesic from $\alpha_i. x_0$ to $\beta_i. x_0$. Then, let $ \textsf{c}':[0, \ell']\f X$ be  the curve 
$$ \textsf{c}' =  \textsf{c}_1'*\cdots * \textsf{c}_m' -    \; \alpha_1.x_0$$
that is,  the concatenation of ($\mathcal Z$-translated of) the geodesics  $ \textsf{c}_i$ such that the endpoint of $ \textsf{c}_1'*\cdots * \textsf{c}_{i-1}'$ 
coincides with the origin of $ \textsf{c}_{i}'$, and  $ \textsf{c}' (0)= x_0$. Finally, let $\gamma'. x_0= [\sum_1^m (\beta_i-\alpha_i)]. x_0$ be the endpoint of $ \textsf{c}'$.  \\
Notice that, as the $ \textsf{c}_i'$ are geodesics, we have:
\small
\vspace{-3mm}

$$ d(x_0, \gamma'.x_0) \leq \sum_{i=1}^m \ell ( \textsf{c}_i') 
\le 2n \diam ( \mathcal Z \backslash X) +  \sum_{i=1}^m d( \textsf{c}(a_i),  \textsf{c}(b_i)) 
\le   2n \cdot  \diam  ( \mathcal Z \backslash X) + \frac{\ell}{2}$$

\normalsize
\noindent Moreover, we have  
 $F(\gamma'. x_0)=\sum_{i=1}^m (F(\beta_i. x_0)- F(\alpha_i . x_0))$ as $F= f$ on $\mathcal Z x_0$,  so
\small
$$\left\| F(\gamma'.x_0)- F(\gamma. x_0) \right\|_{euc}
 =\left\| \sum_i \left[ F(\beta_i.x_0)- F(\alpha_i. x_0) \right]- \sum_i \left[  \textsf{c}_o(b_i)- \textsf{c}_o(a_i) \right] \right\|_{euc} \le$$
 
 \vspace{-4mm}
 
$$\le\sum_{i=1}^m \left\|  F(\beta_i. x_0)- F( \textsf{c}(b_i)) \right\|_{euc} +\sum_{i=1}^m \left\| F(\alpha_i.x_0)- F( \textsf{c}(a_i)) \right\|_{euc} 
 \le 2n M''\cdot \diam  ( \mathcal Z \backslash X)$$
\normalsize
and from this and the Lipschitz property of $f^{-1}$ we deduce that: 
\small
$$d(\gamma'. x_0, \gamma. x_0)=d(f^{-1}(F(\gamma'.x_0)), f^{-1}(F(\gamma. x_0)))
\le  2nM'M''  \cdot \diam  ( \mathcal Z \backslash X) .$$
\normalsize
Then, $d(x_0, \gamma. x_0)\le d(x_0, \gamma'.x_0)+d(\gamma'.x_0,\gamma.x_0) 
\le    \frac{\ell}{2} + 2n \left(  M'M'' +1 \right) \cdot \diam  ( \mathcal Z \backslash X)  $
\normalsize 
\noindent that is, $d(x_0, \gamma. x_0)\le \frac12 d(x_0, 2\gamma. x_0) +  M'''$
 for a constant $M'''=M'''( n,D,\Omega,\sigma)=2n(M'M'' +1) \diam (\mathcal Z \backslash X)$ which is given explicitly by  (\ref{eqM}), (\ref{eqM'}) and Lemma \ref{est_quot}.\\
This implies  the announced inequality 
$\Big|\, d(x_0,\gamma .x_0)- \|\,\gamma\,\|_{st} \, \Big|\le M'''$
for all $\gamma\in\mathcal Z$. 
To get the inequality for all $\gamma\in\Gamma$, let $\gamma_0.x_0$ be a point of $\mathcal Z .x_0$ closest to $\gamma. x_0$; then \linebreak
 \vspace{-9mm}

$$\Big|\, d(x_0,\gamma .x_0)-\|\,\gamma\,\|_{st} \, \Big| \le \Big| d(x_0,\gamma_0.x_0)- \|\,\gamma_0\,\|_{st}\Big|+ 2\,d(\gamma.x_0,\gamma_0. x_0)\le c(n,D,\Omega,\sigma)$$
 for $c(n,D,\Omega,\sigma) = 2 \diam ( \mathcal Z \backslash X) ( nM'M'' +n +1)$.$\Box$
 \vspace{2mm}
 
\noindent In the next section we show that the constant $c$ actually does not depend on $\sigma$.

\vspace{3mm}
\section{Stable systole and asymptotic volume}

We prove here the two relations of (almost) inverse proportionality between $\omega(\Gamma, d)$ and $\stsys(\Gamma, d)$. First notice that, as $\|\,\cdot\,\|_{st}$ is a true norm,  the ball of radius $2D$ in $(\Gamma, \|\,\cdot\,\|_{st})$ is compact; so, there exists $\gamma_1\in\Gamma=\Z^n$  realizing the stable systole.
Let $\sigma=\|\,\gamma_1\,\|_{st}=\stsys(\Gamma, d)$ and  $D=\codiam(\Gamma, d) = \diam (\Gamma \backslash X)$. 

\subsection{Proof of the Abelian Margulis Lemma, upper bound.} 
\label{secmargulisupper}
${}$

\noindent Let ${\mathcal D}_{st}=\{ p \; | \; \| p \|_{st} < \| p- \gamma .o\|_{st} \,\} $ and $\widehat{{\mathcal D}}_{st}=\{ p \; | \; \| p \|_{st} \leq \| p- \gamma .o\|_{st} \,\} $ be respectively the  {\em open} and {\em closed} Dirichlet domains of $\Gamma$ acting on $\R^n$, centered at the origin, with respect to the stable norm, and let $M \geq \diam (\widehat{{\mathcal D}}_{st})$. 
Notice that, in general, the closure of $ {\mathcal D}_{st}$ might be strictly included in 
$\widehat{{\mathcal D}}_{st}$, and    neither $\overline{ {\mathcal D}}_{st}$ nor  $\widehat{{\mathcal D}}_{st}$ a priori tile $\R^n$ under the action of $\Gamma$ (think for instance to the Dirichlet domain of $2 \Z \times \Z$ acting on $(\R^2, \| \; \|_{\infty})$). So, let ${\mathcal F}$ be a closed fundamental domain such that ${\mathcal D}_{st} \subset  {\mathcal F} \subset \widehat{{\mathcal D}}_{st} $; that is,   $\bigcup_{\gamma\in\Gamma}\gamma. {\mathcal F} = \R^n$ and 
$\gamma.\mathring{{\mathcal F}}\cap \gamma'. \mathring {\mathcal F}=\varnothing$ for   $\gamma\neq\gamma'$.
The open ball $B_{st} (r)=\{ p \; | \; \| p \|_{st} < r\} $ is included in $ {\mathcal D}_{st}$ for $r=\frac{\sigma}{2}$, 
so
\vspace{-3mm}

\small

\begin{eqnarray}
\label{equpper}
\frac{\# \left[ B_{(X,d)} (x_0, R) \cap \Gamma \right]}{R^n}  
\leq   \frac{\# \left[ B_{st} (R) \cap \Z^n \right]}{R^n}
\leq  \frac{\Vol(B_{st}(R+M))}{\Vol({\mathcal F}) \cdot R^n}  \leq \hspace{20mm}
\\
\label{equpper2}  \hspace{30mm} \leq   \frac{\Vol ( B_{st} (R+M))}{ \Vol ({\mathcal D}_{st}) \cdot R^n }
 \le  \frac{\Vol ( B_{st} (R+M)) }{ \Vol  ( B_{st} (\frac{\sigma}{2})) \cdot R^n} = \frac{2^n(R+M)^n}{ \sigma^n \cdot R^n}  
\end{eqnarray}
\normalsize
 and taking limits for $ R \rightarrow \infty$ yields the announced inequality.\\ 
Clearly, this inequality is an equality for the standard lattice  $\Z^n$  in  $(\R^n,\| \; \|_\infty)$, but this is not the only case in which the equality is satisfied. Actually, assume  that the equality $\omega(\Gamma,d)=\frac{2^n}{\sigma^n}$ holds: then, all the inequalities in (\ref{equpper}) and (\ref{equpper2})  are equalities for $R\rightarrow \infty$, so 
$\Vol  ( B_{st} ( \frac{\sigma}{2})) = \Vol ({\mathcal D}_{st}) = \Vol ({\mathcal F})$. 
Since $B_{st} ( \frac{\sigma}{2})  \subset  {\mathcal D}_{st} \subset  \mathring {\mathcal F} $, we deduce that  
$B_{st} (\frac{\sigma}{2}) = {\mathcal D}_{st} = \mathring {\mathcal F}$.
This implies that $\overline{{\mathcal D}}_{st} \subset {\mathcal F}$ is a convex  set (being the closure of a ball) which tiles $\R^n$ under the action of $\Gamma$. Actually, assume that there exists
$ p \in \R^n \setminus \Gamma    .\overline{{\mathcal D}}_{st}$. 
Then,  $\R^n \setminus  \bigcup_{ \|   \gamma   \|_{st} \leq  \| p\| + 2M}  \gamma   .\overline{{\mathcal D}}_{st} $ is a non-empty open set, containing a small ball $B_{st} (p, \varepsilon)$ centered at $p$. As ${\mathcal F}$ tiles, there exists $\gamma$ such that 
$\Vol ({\mathcal F} \cap \gamma B_{st}(p,\varepsilon)) \neq 0$. 
This yields a contradiction, as $\overline{{\mathcal D}}_{st} \subset  {\mathcal F} \setminus \gamma B_{st}(p,\varepsilon)$ but 
$\Vol (\overline{{\mathcal D}}_{st}) = \Vol ({\mathcal F})$.

\noindent We show now that $\overline{{\mathcal D}}_{st}$ is a polyhedron. For this, let us first show that  the topological boundary 
$\partial  {\mathcal D}_{st}$ is covered by a finite number of hyperplanes: actually, as the closed sets $\gamma \overline{ {\mathcal D}}_{st}$ tile, we have 
$$\partial  {\mathcal D}_{st}
 = \bigcup_{0<\| \gamma \| \leq 2M} \left( \partial  {\mathcal D}_{st} \cap  \gamma . \partial{\mathcal D}_{st} \right)
    =   \bigcup_{0<\| \gamma \| \leq 2M}  \left(  \overline{ {\mathcal D}}_{st} \cap \gamma. \overline{ {\mathcal D}}_{st}  \right)$$
and as  $\overline{ {\mathcal D}}_{st} \cap \gamma. \overline{ {\mathcal D}}_{st} $ is a convex set with zero measure, it is contained in an affine hyperplane   $H_\gamma = \{ p \, | \, f_\gamma (p) = 1 \}$, for some linear function $f_{\gamma}$; since $\overline{ {\mathcal D}}_{st} $ is convex,  we may assume that   $ \overline{ {\mathcal D}}_{st} \subset H_{\gamma}^{-}$, where  $H_{\gamma}^{-}$ denotes the sub-level set $f_{\gamma} \leq 1$. \linebreak 
 Let $\Gamma_0$ be the subset of nontrivial elements  $\gamma \in \Gamma$ such that $\overline{ {\mathcal D}}_{st} \cap \gamma. \overline{ {\mathcal D}}_{st} \neq \emptyset$.\linebreak 
It then follows that  $\overline{{\mathcal D}}_{st} = \bigcap_{\gamma \in \Gamma_0} H_{\gamma}^-$.
The inclusion  $\overline{{\mathcal D}}_{st} \subset \bigcap_{\gamma \in \Gamma_0} H_{\gamma}^-$ is clear.\linebreak 
On the other hand, given $p \in  \bigcap_{\gamma \in \Gamma_0}  H_{\gamma}^-$,  if $p \not\in  \overline{{\mathcal D}}_{st}$ then the segment  $\overline{op}$ intersects $\partial{{\mathcal D}}_{st}$ at some point $tp$, for   $0<t<1$, hence there exists some $f_\gamma$ such that $f_\gamma (tp) =1$; hence  $f_\gamma (p) > 1$, a contradiction.  
 This shows  that $\overline{{\mathcal D}}_{st}=\overline{ B_{st} ( \frac{\sigma}{2}) } $  is a convex polyhedron tiling $\R^n$  under the action of $\Gamma$, i.e.    a $\Gamma$-parallelohedron (and $B_{st} (1)$ as well). \\
Finally, as  $\big| d - \|\;  \|_{st} \big| < c(n,D,\Omega)$ 
on  $\Gamma.x_{0}$ by the QBD Theorem, by identifying the orbit $\Gamma .x_0$ with  $\Z^n$ we deduce a $\Gamma$-equivariant map $f: (X,d) \rightarrow (R^n,\| \;\|_{st})$ which is a $C$-almost isometry for 
$C=c(n,D,\Omega) + 2 \codiam (\Gamma,d) + 
  \codiam (\Z^n,  \|\cdot\|_{st} ) \le c(n,D,\Omega) + 2D + \sigma$.$\Box$\\

\vspace{-2mm}
\subsection{Proof of the  Abelian Margulis Lemma, lower bound.} 
${}$

\noindent As $\gamma_1$ realizes the stable systole,  for any $\varepsilon>0$ there exists a  $K_\varepsilon$ such that 
$$(1-\varepsilon)\,|k|\,\sigma \le  | \gamma_1^k|_{x_0 }\le (1+\varepsilon)\,|k|\,\sigma  \mbox{  for all } |k|> K_\varepsilon.$$
Complete $\gamma_1$ to a set $\Sigma_n=\{\gamma_1,\gamma_2,...,\gamma_n\}$ of $n$ linearly independent vectors, taking $\gamma_2,...,\gamma_n$ from the generating set $\Sigma'_D \!\!=\!\! \{\gamma\in\Gamma\,|\, d(x_0, \gamma. x_0)\le 2D\}$, and let $\mathcal Z\!=\!\langle\Sigma_n\rangle$. \linebreak
Then, consider the norm   $\|  \;\;\|_{\sigma,2D}$ given by the weighted $\ell_1$-norm on $\R^n$,  relative to the basis $\Sigma_n$, with   weights $\ell(\gamma_1)=\sigma$ and $\ell(\gamma_i)=2D$ for $i\neq1$.  
Finally, let $ \mathcal Z_{\varepsilon}:=\{ \gamma_1^{k_1}\cdots\gamma_n^{k_n}\,|\, |j_1|> K_\varepsilon\}$. 
Then, for all $\gamma \in  \mathcal Z_{\varepsilon}$ we   have:
\vspace{-5mm}

%


$$| \gamma |_{x_0} 
\leq  | \gamma_1^{k_1} |_{x_0} + \sum_{k=2}^n   | \gamma_1^{k_i} |_{x_0}
\leq (1+\epsilon) |k_1| \cdot \sigma  + 2D  \sum_{k=2}^n   | k_i |  
\leq   (1+\epsilon)   \cdot \|  \gamma \| _{\sigma,2D}.$$

\noindent Therefore we obtain:
\vspace{-3mm}
\small


 \begin{equation}\label{eqlower}
\omega(\Gamma, d_{x_0})
\ge \omega(\mathcal Z, d_{x_0}|_{\mathcal Z})
\ge \omega(\mathcal Z_{\varepsilon}, d_{x_0}|_{\mathcal Z_{\varepsilon}}) 
 \ge \frac{ \omega(\mathcal Z_\varepsilon,  \| \;\;\| _{\sigma,2D}  )}{(1+\varepsilon)^n  } 
  =\frac{  \omega(\mathcal Z_\varepsilon, d_{\Sigma_n})}{(1+\varepsilon)^n\,\sigma(2D)^{n-1}}
  \end{equation}

\normalsize
\noindent which gives the announced bound, as $\epsilon>0$ is arbitrary  and since 
\vspace{-2mm}
\small

$$\omega(\mathcal Z_\varepsilon, d_{\Sigma_n}) \! =\! \omega(\mathcal Z, d_{\Sigma_n})\! =\! \omega(\Z^n, \| \;\ \|_1) \! = \!  \frac{2^n}{n!} \; ;$$ 
\normalsize

\noindent actually, 
 the set $\mathcal Z \setminus \mathcal Z_\varepsilon$ has polynomial growth of order $n-1$ and is negligible in the computation of the asymptotic volume, while $(\mathcal Z, d_{\Sigma_n})$ is isometric to $\Z^n$ \linebreak with the canonical word metric.  \\
Notice that the equality holds for the action  the standard lattice $\Z^n$ on $(\R^n, \|\,\,\|_1)$.\linebreak 
Best, assume that, for $\Gamma$ acting on $(X, d)$, we have the equality
$\omega(\Gamma, d_{x_0})= \frac{2}{n!   D^{n-1}\sigma}$. 
In particular, the first inequality in (\ref{eqlower})   is an equality, which implies  $[\Gamma: \mathcal Z] = 1$.
Moreover, we deduce that  $\omega(\Gamma, \| \;\; \|_{st} ) =\omega(\Gamma, d_{x_0})= 
 \omega(  \Gamma,  \| \;\;\| _{\sigma,2D}) $.
However,   by construction, the stable and weighted norms satisfy $  \| \;\;\| _{\sigma,2D} \geq \| \; \;\|_{st}$; then, being norms, we know that the equality of asymptotic volumes implies the equality of $1$-balls, so $  \| \;\;\| _{\sigma,2D} = \| \; \;\|_{st}$. Therefore $\| \; \;\|_{st}$ is  affine equivalent to the $\ell_1$ norm  $\|\cdot\|_1$, via an affine map sending $\Z^n$ to the lattice $\Gamma_0=\sigma \!\cdot\! \Z \times 2D \!\cdot \! \Z^{n-1}$ of $\R^n$.
It follows by the QBD Theorem that the action of $\Gamma$ on $(X,d)$ is equivalent, via an equivariant $C$-almost isometry 
$ f: (X,d) \rightarrow (\R^n, \|\cdot\|_1)$, to the action of $\Gamma_0$ on $(\R^n, \|\cdot\|_1)$, for 
$C=c(n,D,\Omega) + 2 \codiam (\Gamma,d)  +  \codiam (\Gamma_0,  \|\cdot\|_1)
\le c(n,D,\Omega) + (2n+2)D$.$\Box$
 \vspace{2mm}

\begin{rmk}\label{remstimec}
Let $\sigma = \sys (\Gamma, d) \geq \sigma$ and $ \omega(\Gamma, d) \leq \Omega$ as in section \S2.
\vspace{1mm}

\noindent (i) Using the lower bound  given by  the Abelian Margulis Lemma and the fact that $D \geq \frac{\sigma}{2}$,
we find $\Omega D^n \geq \omega(\Gamma, d) D^n\geq  \frac{1}{n!}$. This estimate, together with (iii) of Lemma \ref{est_quot},   plugged in the expressions   (\ref{eqM}),   (\ref{eqM'}) for $M, M'$, and in the expressions for $M'', M'''$ and 
 $c= M'''+ 2 \diam  ( \mathcal Z \backslash X)$ of \S\ref{sectionproof},  yields  the following estimate  for the   constant $c$ of the QBD Theorem:
 \vspace{-4mm}

$$  c (n,D, \Omega, \sigma) = c(n,D, \Omega) \leq 2^{n^2 + 6n + 10} \cdot n^2  \cdot (n! )^{n+2} \cdot D ( \Omega D^n +1)^{n+4}$$

\noindent  Notice that the quantity $\Omega D^n$ is scale invariant. 
\vspace{1mm}

\noindent (ii) We also remark, for future reference, that the same  computations show that the constant $c$ that we find is   $\gg nD$, namely
$ c(n,D, \Omega) \geq  2^{n^2 + 6n + 8} n^2  (n! )^n D$.

\end{rmk}

\vspace{1mm}
\begin{rmk}
As a consequence, we have the explicit bound
\small
$$\left| \frac{  | \gamma |_{x_0} }{ \|\, \gamma \,\|_{st}  } -1 \right|
\! \leq \! \frac{c(n,D,\Omega)}{ \|\, \gamma \,\|_{st}  } $$
\normalsize
This should be compared with an asymptotics given by  Gromov in \cite{gro2}:
\small
\begin{equation}
\label{eqgromov}
\left| \frac{  | \gamma |_{x_0} }{ | \gamma |_{H_1}  } -1 \right|
\leq \frac{c_{\bar X}}{ | \gamma |^{n-1}_{H_1}}
\end{equation}
\normalsize
for the mass  of $\gamma \in H_1 (\bar X, \Z)$.
Notice however that Gromov's bound is purely qualitative (no information can be deduced on the constant  $c_{\bar X}$ from his argument)  and that we always have $\|\, \gamma \,\|_{st} \le | \gamma |_{H_1} $, by the characterization of the stable norm in real homology recalled in the introduction. 
\end{rmk}

\vspace{3mm}
\section{Examples}

Here we show that  the  constant $c=c(n,D,\Omega)$ of Theorem \ref{QBDT} necessarily depends on each of the three parameters rank,  diameter and asymptotic volume. \linebreak
\noindent We say that a sequence of actions of torsionless, discrete abelian groups $\Gamma_k$ on $(X_k,d_k)$ is {\em noncollapsing} if there exists $\sigma >0$ such that $\stsys( \Gamma_k, d_k ) > \sigma$ for all $k$.

\begin{exmp} {\em Collapsing actions with fixed rank and bounded co-diameter.} 
\label{exdiameter}

 \noindent Let $\Z$ act on $(X_k, d_k) = \mathcal{C} (\Z, S_k)$, the Cayley graph of $\Z$ with respect to the generating set  $S_k= \{ \pm1, \pm k \}$, and let $\| \; \|_{st,k}$ be the associated stable norm. Then:
\vspace{-1mm}

  \noindent (i) $\codiam( \Z, d_k ) = 1$;
    
\noindent  (ii) $\sys (\Z, d_k) = 1$, while $\stsys (\Z, d_k) \stackrel{k \rightarrow \infty}{\longrightarrow} 0$, as
 $|1|_{st,k} = \lim_{m \rightarrow \infty} \frac{d_k (0,km)}{km} \leq \frac{1}{k}$.
 \vspace{1mm}
  
\noindent  (iii) $\omega (\Z, d_k) \rightarrow \infty$, as a consequence of Lemma  \ref{QBDT};
 \vspace{1mm}

\noindent  (iv) $d_{2k} (0,k) = k$, while $\|k\|_{st, 2k} \leq \frac12$.
 \vspace{1mm}

\noindent This example shows that, the rank and the co-diameter of $(\Gamma_k, d_k)$ being fixed, without any assumption on the asymptotic volume (or  the stable systole) the difference between the distance and the associated stable norm can be arbitrarily large. 
\\
It also shows that, whereas the collapse of the systole forces the asymptotic  volume to diverge (by the Abelian Margulis Lemma), the converse is not true. \\
Notice that, with little effort,  the example can be modified into a sequence of  $\Z$-coverings  of a compact Riemannian manifold   with the same properties, in the following way.  
Start with the $\epsilon$-tubular neighbourhood in $\R^3$ of a bouquet of two circles $\alpha, \beta$ with length $1$, and  consider its boundary $\bar Y$. 
Let $(Y_k, d_k)$ the  Riemannian covering of $\bar Y$ associated to the subgroup $N= \langle \alpha, \alpha^k \beta^{-1}\rangle$ of $H_1 (\bar Y, \Z)$: then, there exists a $(1+\delta, \delta)$-quasi isometry between $Y_k$ and the above graph $X_k$, with $\delta \approx \epsilon$, which is equivariant with respect to the actions of $\Gamma = H_1 (\bar Y, \Z)/N \cong \Z$. Therefore, $\Gamma$ acts on $Y_k$ with the same properties (up to multiplicative costants $1+\delta$ in the above estimates (i)-(ii)-(iv)).


\end{exmp}

\begin{exmp} {\em Noncollapsing actions with fixed rank and large co-diameter.} 
\label{exdiameter}

 \noindent Let $\Z$ act on $(X_k, d_k) = k \cdot \mathcal{C} (\Z, S_p)$, the Cayley graph of $\Z$ with respect to the generating set 
 $S_p=\{ \pm 1, \pm p \}$, with $p>1$ fixed, and the graph metric  dilatated by a factor $k$. Let $\| \; \|_{st,k}$ be the associated stable norm.  Then:
  \vspace{1mm}
  
  \noindent (i) $\diam( X_k, d_k ) = k$;
    \vspace{1mm}
    
\noindent  (ii) $\stsys (\Z, d_k) \geq \frac{k}{p}$, since 
$\|m\|_{st,k} = \lim_{h\rightarrow \infty} \frac{d_k (0,mhp)}{hp}
 = k \cdot\frac{m}{p} \geq \frac{k}{p}$.
 \vspace{1mm}
  
\noindent  (iii) $ \frac{p}{k} \leq \omega (\Z, d_k) \leq\frac{2p}{k} $, as a consequence of Lemma  \ref{margulis};
 \vspace{1mm}
  
\noindent  (iv) $d_{k} (0,1) = k$, while $\|1\|_{st, k} =\frac{k}{p}$.
 \vspace{1mm}

\noindent This example shows that, the rank and the asymptotic volume being bounded, without any assumption on the diameter the difference between the distance and the associated stable norm can be arbitrarily large. 
 
\end{exmp}

\begin{exmp} {\em Noncollapsing actions with  large rank and bounded co-diameter.} 
\label{exrank}

 
\noindent Consider a  round sphere $(S^3,d)$  with north pole $x_0$,  and remove  an arbitrarily large number  $n$ of  small, disjoint balls $B_i$,  centered at $m$ equatorial points, 
with boundary $2$-spheres $S_i$  
(so that $d (x_0,S_i) \sim \frac{\pi}{2}$); 
 then, take $n$ copies  $T_i$ of a flat torus,  each with a small ball $B'_i$ removed and  boundary spheres $S'_i$,  glue   the (almost isometric) spheres  $S_i, S'_i$ through a cylinder of length $\ell$, and smooth the metric to obtain a Riemannian manifold $\bar X_n$. We may assume that $\ell$ is much larger than 
 the length $\sigma$ of the shortest nontrivial 1-cycle  in the flat  torus $T_i$ 
 (which realizes the stable systole of $H_1(T_1, \Z)$ acting on  the universal covering of $T_i$),
  and that, nevertheless,  $\diam (\bar X_n)$ stays bounded.
The groups $\Gamma_n \! = \! H_1 (\bar X_n, \Z) \! = \! \bigoplus_{i=1}^n  H_1 ( T_i, \Z) \! \cong\!  \Z^{3n}$ then act  on the Riemannian homology coverings $(  X_n,  d_n)$ of $\bar X_n$ without collapsing: actually, 
 any class $\gamma_i \in H_1 (T_i, \Z)$ has   length $\ell (\gamma_i)$ in  $X_n$  not smaller than its original length in $T_i$ 
 (the ball  $B'_i$ has been replaced by an almost  flat cylinder), and  any decomposable class $\gamma = \sum_i \gamma_i$ with $\gamma_i \in H_1 (T_i, \Z)$  has length greater than $\sum_i \ell (\gamma_i)$.
Thus,   $\stsys (\Gamma_n, d_n) \geq \sigma$ for all $n$. 
On the other hand,  for every $\gamma= \sum_i \gamma_i $, with  nontrivial  components $\gamma_i \in H_1 (T_i, \Z)$ for all $i$, we have 
$d(x_0, \gamma x_0) \geq  2 n \ell + \sum \ell (\gamma_i)$ (as the shortest geodesic loop representing $\gamma$ must travel forth and back at least $n$ cylinders), 
while $ \| \gamma \|_{st} \leq \sum_i  \ell ( \gamma_i)$;
hence $ d(x_0, \gamma x_0) - \| \gamma \|_{st}$ diverges for $n \rightarrow \infty$.
Notice that in these examples  the asymptotic volume  $\omega (\Gamma_n,d_n)$ stays bounded for $n \rightarrow \infty$, by Lemma \ref{margulis}, although the rank is arbitrarily large.
   
\end{exmp}

\vspace{3mm}
We conclude this section with an example showing that the Bounded Distance Theorem may fail for abelian actions on metric spaces which are not length spaces. 
{\em Inner metric spaces}
\footnote{Any finitely generated group $\Gamma$ endowed with a word length, or with a  geometric distance deduced from a cocompact action on a length space, is an inner metric space.},  as defined by P. Pansu \cite{pansu}, are the closest spaces to  length spaces:  $(X,d)$ is inner if, for every $\epsilon>0$, there exist $\ell (\epsilon)$ such that for all $x,x'\in X$ there exists a sequence of points $x_0=x,x_1,\cdots, x_{N+1}=x'$ with $d(x_i,x_{i+1})\le\ \ell(\epsilon)$ and $\sum_{i=1}^{N+1} d(x_{i-1},x_i) \leq (1+\epsilon) d(x,x')$. 
 The following is the simplest example of inner metric space
where the Bounded Distance Theorem does not hold:

\begin{exmp}{\em Noncollapsing $\Z$-actions on inner spaces with bounded co-diameter.}
\label{sqrt-d}

\noindent Consider the group $\Z$ endowed with the left invariant metric induced by the norm 
$|\!|\!| m |\!|\!|=|m|+\sqrt{|m|}$.
 It is straightforward to check that $|\!|\!| \;  |\!|\!|$ defines an inner metric  on $\Z$. Actually, given $\epsilon>0$, choose an integer $\ell > 4/\epsilon^2$, and write any $m\in \N$ as  $m=N\ell+r$, with $r<\ell$. If $N=0$, there is nothing to prove; otherwise call $x_i = i\ell$ for $i\leq N$ and $x_{N+1}=m$, so
 \vspace{-2mm} 
  
 \small
 
$$\hspace{-2cm}
\frac{  \sum_{i=1}^{N+1} |\!|\!| x_i - x_{i-1} |\!|\!|   }{ |\!|\!| m  |\!|\!|}
\leq \frac{  \sum_{i=1}^{N} |\!|\!| \ell |\!|\!|  + |\!|\!| r |\!|\!|}
          {  |\!|\!| N\ell+r |\!|\!|  }
=   \frac{  N(\ell+\sqrt{\ell})+ r+\sqrt{r}   }{   (N\ell+r)+\sqrt{N\ell+r}  } \le$$

\vspace{-1mm}

$$\hspace{2cm}\le \frac{  N\ell+r  }{ (N\ell+r)+\sqrt{N\ell+r}  }
    +\frac{1}{\sqrt{\ell}}+ \frac{\sqrt{r}}{N\ell}
\le 1+\frac{2}{\sqrt{\ell}}<1+\epsilon$$

\normalsize
\noindent Then, $\Z$ acts by left translation on itself, and 
the stable norm associated to $|\!|\!| \;  |\!|\!|$ coincides with the absolute value $|\,|$.
Therefore, we have 
$\codiam (\Z, |\!|\!| \; |\!|\!|)=\frac12 + \sqrt{\frac12}$ and
$\stsys(\Z, |\!|\!| \; |\!|\!|) =1$, but 
$|\!|\!| m |\!|\!| - |\!|\!| m |\!|\!|_{st}=\sqrt{|m|}$
is not bounded.

\end{exmp}

\vspace{1mm}
\section{On the number of connect components of optimal cycles}
\label{sectioncc}

Let $\bar X$ be a Riemannian manifold, with torsion free  homology covering $(X, d)$.
Let  $x_0 \in X$ be fixed, let   $d_{x_0}$ be the induced distance on  $\Gamma = H_1 (\bar X, \Z)$  acting on $(X,d)$,
and $\| \; \|_{st}$ be the associated stable norm on $H_1 (\bar X, \R)$ as explained in \S1. \linebreak
Let $\lambda$ be the Busemann  measure of the normed space $\left( H_1 (\bar X, \R), \|\cdot\|_{st} \right)$, that is the Lebesgue measure  assigning to its unit ball $B_{st} (1)$    the volume of the unitary euclidean $n$-ball (which coincides with the $n$-dimensional Hausdorff measure).

\noindent Finally, let ${\mathcal F}$ be a closed fundamental domain  included in the closed Dirichlet domain $\widehat{ {\mathcal D}}_{st}$ centered at the origin, for    $\Gamma$  acting on $\left(H_1 (\bar X, \R), \|\cdot\|_{st} \right)$,  as in \S\ref{secmargulisupper}.


 \vspace{1mm}
 As recalled in the introduction,  an easy packing of fundamental domains shows that  the asymptotic volume of the measure metric space $(H_1 (\bar X, \R), \| \; \; \|_{st}, \lambda)$ is 
\vspace{-3mm}

\begin{equation}
\label{omegalambda} \omega_\lambda (\R^n, \|\;\;\|_{st})  
    = \omega (\Gamma, \|\;\|_{st}) \cdot \lambda (\mathcal F) 
 \end{equation}

\noindent and, since  $\|\;\;\|_{st}$ is a norm,  this also equals 
the volume $\lambda (B_{st} (1))$ of the unit ball.

\noindent Then, as a consequence of the Bounded Distance Theorem (even without any estimate of the constant $c$),
one gets: 
%
$\omega (\Gamma, d) =  \omega (\Gamma, \| \; \|_{st}) =   \lambda ( B_{st}(1) )  / \lambda (\mathcal F)   $.\linebreak
%
\noindent Let  us call   $V=\lambda (\mathcal F)$ and $V_1=\lambda \left(B_{st}(1)\right)$, so that $\omega (\Gamma, d) = V_1/V$.

 \vspace{1mm}
  Let us now fix some notations for balls and annuli  and for the corresponding growth functions.   
  We will write
  \vspace{-3mm}

$$ B_{(\Gamma, d_{x_0})}  (R)
=\left\{ \gamma \in\Gamma  \;|\;    | \gamma |_{x_0} < R  \right\}$$

\vspace{-4mm}

$$ \;\;\;\;A_{(\Gamma, d_{x_0})}  (r, R)
=\left\{ \gamma \in\Gamma  \;|\;  r \le | \gamma |_{x_0} < R  \right\}$$

\noindent and similarly
$B_{(\Gamma, |  \;|_{H_1})}$,  $B_{(\Gamma, st)}$, 
and $A_{(\Gamma, |  \;|_{H_1})}$,   $ A_{(\Gamma, st)}$ 
for   balls and  annuli    in $\Gamma$ with, respectively,  the mass and the stable norm.
We will  write $v_{_{^{\bullet}}}  (R)$, $v_{_{^{\bullet}}} (r,R)$ for  the corresponding  cardinalities.
Finaly, we will use  $B_{st}$ and $A_{st}$ for ball and annuli in   $(H_1 (\bar X, \R), \|\;\;\|_{st})$, 
and write  $v_{st} (R) = \lambda (B_{st} (R) )$, $v_{st} (r, R) = \lambda (A_{st} (r, R) )$.

 \vspace{2mm}
 A by-product of the QBD Theorem is the following  explicit estimate of  the growth function of annuli in $H_1(\bar X, \Z)$ w.r. to the mass:

\begin{prop}\label{v-estimates}
Assume  $n \! = \! {\rm rank} \, H_1(\bar X, \Z)$,  $\diam (\bar X) \! \leq\!  D$ and $\omega( X,d) \! \leq\!  \Omega$. 

\noindent Let  $c=c(n, D, \Omega)$ be as in the QBD Theorem \ref{QBDT}.
If   $\Delta > 4nD +  c$   we have:
\vspace{-2mm}
 \small
 
\[   n \omega(\Gamma,d)  \cdot  \Delta (R - \Delta)^{n-1}
\le v_{(\Gamma, | \; |_{H_1})} (R-\Delta,R+\Delta )
\le  3 n \omega(\Gamma,d)  \cdot   \Delta (R- 3 \Delta)^{n-1}
\mbox{ for } R\geq 0 \]

\vspace{-3mm}

$$A   (k\Delta)^{n-1}
\le v_{(\Gamma, | \; |_{H_1})} (k\Delta,(k+1)\Delta )
\le  B (k\Delta)^{n-1}
\mbox{ for all } k\geq 1   $$

\normalsize
\vspace{1mm}
\noindent for constants  $A= n  \cdot  \omega(\Gamma,d) \cdot \Delta$ and  $B=3^n \cdot A$.
\end{prop}

\textbf{Proof.} By Theorem \ref{QBDT}  we have
$\|\,\gamma\,\|_{st}\le |\,\gamma\,|_{H_1}\le  |  \, \gamma \,  |_{x_0} \le \| \, \gamma \, \|_{st}+ c$
and thus  
$$A_{(\Gamma, st) } (r, R-c)
\subseteq A_{(\Gamma, | \; |_{H_1})} (r, R)
\subseteq A_{(\Gamma, st) }(r-c,R)$$
Notice  that, if   $D_{st} = \diam_{st} (\mathcal F)$  is  the diameter of $\mathcal F$   with respect to the stable norm,  we have  
  $D_{st} \leq 2nD$. Actually, choose $n$ linearly independent vectors $v_1,...,v_n\in\Sigma'_D$ from   the generating set  
$\Sigma'_D = \{ \gamma \in \Gamma \;: \; |\gamma|_{x_0} \leq 2D\} $: the n-parallelotope ${\mathcal P}$ determined by these vectors clearly  contains   $\mathcal F$, so
$D_{st}   \le \diam_{st}({\mathcal P})
\leq 2n\,D.$ \linebreak 
 By an argument  of packing and covering of the annuli 
 $A_{st}(r, R)$ with copies of  $\mathcal F$,    we obtain
\vspace{-3mm}

 \small
$$\frac{\lambda \left(A_{st}(r+ 2n\,D, R-C- 2n\,D)\right)}{V}
\le v_{(\Gamma, | \; |_{H_1})} (r, R)
\le \frac{\lambda \left(A_{st}(r-2n\,D-C, R+2n\,D)\right)}{V}$$
\normalsize

\noindent We estimate  the left-hand  inequality for $R \geq r + 4nD +c$ as
\small
$$\lambda \left(A_{st}(r+ 2n\,D, R-C-2n\,D)\right)= V_1\cdot\left[(R- C-2n\,D)^n-(r+2n\,D)^n\right]$$ 
$$\hspace{54mm}\ge nV_1\cdot    \left( R-r - (4nD+c) \right) \cdot (r+2nD)^{n-1}$$
\normalsize
and, similarly,  the right-hand   as
\small
$$\lambda \left(A_{st}(r-C-2n\,D, R+2n\,D)\right)
\le nV_1\cdot  \left( R-r + (4nD+c) \right) \cdot (R+2nD)^{n-1}.$$
\normalsize
Choosing $\frac{R-r}{2} = \Delta \geq 4nD+c$ we  obtain 
\small
 $$n \frac{V_1}{V} \cdot  \Delta (R - \Delta)^{n-1}
\le v_{(\Gamma, | \; |_{H_1})} (R-\Delta,R+\Delta )
\le  3 n \frac{V_1}{V}  \cdot   \Delta (R- 3 \Delta)^{n-1}$$
\normalsize
which proves  (i). The second statement follows from (i) taking $R=(k+1)\Delta$.$\Box$

\vspace{7mm}
\textbf{Proof of Theorem \ref{teorcc}.}
Let $\Delta$ be as in Proposition \ref{v-estimates} above. \\
 First, we consider the case where $N(\gamma)\le v_{(\Gamma, H_1)} (2\Delta)$.
As $\| \; \: \|_{st} \leq | \; \; |_{H_1}$,  we have  
$$N(\gamma)\le \frac{\lambda \left( B_{st} (2\Delta +D_{st} \right) }{V}
\le \frac{V_1}{V} (2\Delta + nD)^n = \omega (\Gamma, d) (9nD +c)^n$$
and using the explicit estimates for $c$ given in  Remark \ref{remstimec}(i)\&(ii) we get the announced bound $N(n,D, \Omega)$.
Assume now that $N(\gamma)> v_{(\Gamma, H_1)} (2\Delta)$.
Then, there exists  $m=m(\gamma) \ge 1$  such that:
$$\sum_{k=0}^m v_{(\Gamma, H_1)}  (k\Delta, (k+1)\Delta)< N(\gamma)\le\sum_{k=0}^{m+1} v_{(\Gamma, H_1)}  (k\Delta, (k+1)\Delta).$$
\pagebreak

\noindent Then, using the estimates  of Proposition \ref{v-estimates}  we   find 
\vspace{-2mm}

\begin{equation}\label{est1}
N(\gamma) \le  B  \sum_{k=0}^{m+1} (k \Delta)^{n-1} \leq \frac{B\Delta^{n-1}}{n} (m+2)^n.
\end{equation}

\noindent Now, observe that if $ \textsf{c} $ is an optimal cycle representing $\gamma$ with $N(\gamma)$ connected components, then its components $\textsf{c} _i$ are non-homologous to each other; thus, its total length  is at least
$ \ell \geq \sum_{k=0}^{m}  (k \Delta) \cdot v_{(\Gamma, H_1)}  (k\Delta, (k+1)\Delta) $.
Using the estimates  of Proposition \ref{v-estimates}  we   find 
\vspace{-2mm}

\begin{equation}\label{est2}
\ell= |\,\gamma\,|_{H_1} \ge A\sum_{k=0}^{m} (k\Delta)^n  
\ge \frac{A \Delta^n}{n+1}\cdot m^{n+1}.
\end{equation}

\normalsize 
\noindent Putting together the two estimates (\ref{est1}) and (\ref{est2}) above and we obtain:
\small
 \vspace{-3mm}
 
\[ N(\gamma) \le   \frac{B\,\Delta^{n-1}}{n}  (m +2 )^n 
\leq   \frac{B\,\Delta^{n-1}}{n}  \left( \!\! \sqrt[n+1]{ \frac{ (n+1) \ell }{ A \Delta^n} } +2 \right)^n  
\le \frac{3^n  \Delta^{n-1}B }{n \Delta^{\frac{ n^2}{ n+1} }   }\cdot \left( \frac{ (n+1)\ell }{A }  \right)^{n/(n+1)}  
\]

\normalsize 
\noindent  and as   $ A= n \omega(\Gamma, d) \Delta = 3^{-n}B$, this yields  $N(\gamma) 
\le  3^{2n}  \cdot \sqrt[n+1]{(1 + \frac1n)^n  \cdot \omega(\Gamma, d) \ell^n  }. \Box$

\end{document}